\documentclass{crm-article}

\usepackage{amsmath,amsfonts,amssymb,mathtools} 
\usepackage{algorithm}
\usepackage[noend]{algorithmic}
\usepackage{cite}
\usepackage{xcolor}
\usepackage{hyperref}
\hypersetup{
    colorlinks=true,
    linkcolor=blue,
    citecolor=blue,
    urlcolor=blue,
}
\usepackage{xcolor}

\newcommand{\e}{\varepsilon}
\newcommand{\R}{{\mathbb R}} 
\newcommand{\E}{{\mathbb E}}

\newcommand{\la}{\langle}
\newcommand{\ra}{\rangle}

\newcommand{\dom}{\text{dom}}
\newcommand{\N}{\mathcal N}
\newcommand{\hy}{\hat y}
\newcommand{\ty}{\tilde y}
\newcommand{\prox}{\text{prox}}
\newcommand{\A}{\mathcal A}

\begin{document}





\journalVol{10}
\journalNo{1} 
\setcounter{page}{1}

\journalSection{Математические основы и численные методы моделирования}
\journalSectionEn{Mathematical modeling and numerical simulation}

\journalReceived{12.02.2022.}
\journalAccepted{01.06.2016.}

\UDC{519.853.62}
\title{Тензорные методы внутри смешанного оракула для решения задач типа min-min}
\titleeng{Tensor methods inside mixed oracle for min-min problems}
\thanks{Исследование выполнено при поддержке Министерства науки и высшего образования Российской Федерации (госзадание) №075-00337-20-03, номер проекта 0714-2020-0005.}
\thankseng{The research is supported by the Ministry of Science and Higher Education of the Russian Federation (Goszadaniye) №075-00337-20-03, project No. 0714-2020-0005.}

\author[1,2]{\firstname{П.\,А.}~\surname{Остроухов}}
\authorfull{Петр Алексеевич Остроухов}
\authoreng{\firstname{P.\,A.}~\surname{Ostroukhov}}
\authorfulleng{Petr A. Ostroukhov}
\email{ostroukhov@phystech.edu}
\affiliation[1]{Московский физико-технический институт,\protect\\ 
Россия, 141701, Московская область, г. Долгопрудный, Институтский переулок, д.9}
\affiliationeng{Moscow Institute of Physics and Technology,\protect\\ 
9 Institutskiy per., Dolgoprudny, Moscow Region, 141701, Russian Federation}

\affiliation[2]{Институт проблем передачи информации им. А.А. Харкевича Российской академии наук,\protect\\
Россия, 127051, г. Москва, Большой Каретный переулок, д.19 стр. 1}

\affiliationeng{Institute for Information Transmission Problems of Russian Academy of Sciences,\protect\\
Bolshoy Karetny per. 19, build.1, Moscow, 127051 Russian Federation}

\begin{abstract}
    В данной статье рассматривается задача типа min-min: минимизация по двум группам переменных. Данная задача в чём-то похожа на седловую (min-max), однако лишена некоторых сложностей, присущих седловым задачам. Такого рода постановки могут возникать, если в задаче выпуклой оптимизации присутствуют переменные разных размерностей, или если какие-то группы переменных определены на разных множествах. Подобная структурная особенность проблемы дает возможность разбивать её на подзадачи, что позволяет решать всю задачу с помощью различных смешанных оракулов. Ранее в качестве возможных методов для решения внутренней или внешней задач использовались только методы первого порядка или методы типа эллипсоидов. В нашей работе мы рассматриваем данный подход с точки зрения возможности применения алгоритмов высокого порядка (тензорных методов) для решения внутренней подзадачи. Для решения внешней подзадачи мы используем быстрый градиентный метод.
    
    Мы предполагаем, что внешняя подзадача определена на выпуклом компакте, в то время как для внутренней задачи мы отдельно рассматриваем задачу без ограничений и определенную на выпуклом компакте. В связи с тем что тензорные методы по определению используют производные высокого порядка, время на выполнение одной итерации сильно зависит от размерности решаемой проблемы. Поэтому мы накладываем еще одно условие на внутреннюю подзадачу: её размерность не должна превышать 1000. Для возможности использования смешанного оракула нам необходимы некоторые дополнительные предположения. Во-первых, нужно чтобы целевой функционал был выпуклым по совокупности переменных и чтобы его градиент удовлетворял условию Липшица также по совокупности переменных. Во-вторых, нам необходимо чтобы целевой функционал был сильно выпуклый по внутренней переменной и его градиент по внутренней переменной удовлетворял условию Липшица. Также, для применения тензорного метода нам необходимо выполнение условия Липшица $p$-го порядка ($p > 1$). Наконец, мы предполагаем сильную выпуклость целевого функционала по внешней переменной, чтобы иметь возможность использовать быстрый градиентный метод для сильно выпуклых функций.
    
    Стоит отметить, что в качестве метода для решения внутренней подзадачи при отсутствии ограничений мы используем супербыстрый тензорный метод. При решении внутренней подзадачи на компакте используется ускоренный проксимальный тензорный метод для задачи с композитом.
    
    В конце статьи мы также сравниваем теоретические оценки сложности полученных алгоритмов с быстрым градиентным методом, который не учитывает структуру задачи и решает её как обычную задачу выпуклой оптимизации (Замечания \ref{rem:fgm vs our method} и \ref{rem:fgm vs our method on compact}).
\end{abstract}


\begin{abstracteng}
    In this article we consider min-min type of problems or minimization by two groups of variables. In some way it is similar to classic min-max saddle point problem. Although, saddle point problems are usually more difficult in some way. Min-min problems may occur in case if some groups of variables in convex optimization have different dimensions or if these groups have different domains. Such problem structure gives us an ability to split the main task to subproblems, and allows to tackle it with mixed oracles. However existing articles on this topic cover only zeroth and first order oracles, in our work we consider high-order tensor methods to solve inner problem and fast gradient method to solve outer problem.
    
    We assume, that outer problem is constrained to some convex compact set, and for the inner problem we consider both unconstrained case and being constrained to some convex compact set. By definition, tensor methods use high-order derivatives, so the time per single iteration of the method depends a lot on the dimensionality of the problem it solves. Therefore, we suggest, that the dimension of the inner problem variable is not greater than 1000. Additionally, we need some specific assumptions to be able to use mixed oracles. Firstly, we assume, that the objective is convex in both groups of variables and its gradient by both variables is Lipschitz continuous. Secondly, we assume the inner problem is strongly convex and its gradient is Lipschitz continuous. Also, since we are going to use tensor methods for inner problem, we need it to be $p$-th order Lipschitz continuous ($p > 1$). Finally, we assume strong convexity of the outer problem to be able to use fast gradient method for strongly convex functions.
    
    We neet to emphasize, that we use superfast tensor method to tackle inner subproblem in unconstrained case.  And when we solve inner problem on compact set, we use accelerated high-order composite proximal method.
    
    Additionally, in the end of the article we compare the theoretical complexity of obtained methods with regular gradient method, which solves the mentioned problem as regular convex optimization problem and doesn't take into account its structure (Remarks \ref{rem:fgm vs our method} and \ref{rem:fgm vs our method on compact}).
\end{abstracteng}
\keywordeng{tensor methods}
\keywordeng{high-order smoothness}
\keywordeng{strong convexity}
\keywordeng{mixed oracle}
\keywordeng{inexact oracle}

\maketitle

\paragraph{Введение}
    
    На данный момент в оптимизации существует множество методов для различных классических постановок (седловые задачи, выпуклая оптимизация), которые являются наиболее общими. Для этих задач известны нижние оценки \cite{nesterov2004introduction, nemirovsky1983problem}, и известны (суб)оптимальные алгоритмы, достигающие этих нижних оценок. Для дальнейшего ускорения имеющихся алгоритмов научное сообщество всё больше начинает смотреть на структуру имеющейся задачи. 
    Если говорить о выпуклой оптимизации, то сравнительно недавно возникла задача типа min-min, которая схожа с седловой задачей, хотя не обладает некоторыми сложностями, которые возникают при решении седловых задач. Однако, задача с подобной структурой является относительно новой и еще недостаточно изучена \cite{jungers2011min,konur2017set,bolte2020ah,gladin2021solving_a,gladin2021solving_b}. Основная мотивация для решения таких задач заключается в транспортных приложениях \cite{gasnikov2020traffic}.
    Формально стандартная постановка задачи типа min-min выглядит следующим образом:
    \begin{equation}\label{eq:min-min statement}
        \min_{x \in Q_x} \min_{y \in Q_y} F(x, y),
    \end{equation}
    где $Q_x \subseteq R^m$ и $Q_y \subseteq \R^n$ -- некоторые непустые выпуклые множества, $F(x, y)$ является выпуклой по совокупности переменных. 
    
    В работах \cite{gladin2021solving_a,gladin2021solving_b} рассматривается возможность применения смешанного оракула к данной задаче: задача разбивается на внутреннюю и внешнюю. Обе задачи решаются методами первого или нулевого порядка. Но, так как внутренняя задача решается с определенной точностью, в решении внешней задачи использется неточный оракул. В нашей работе также используем концепцию смешанного оракула, однако вместо методов нулевого порядка исследуем применимость методов $p$-го порядка, $p > 1$.
    
    Как известно, методы нулевого порядка имеют более высокую скорость сходимости, однако при большой размерности стоимость одной итерации становится слишком высока по сравнению с градиентными методами. К примеру, у метода эллипсоидов скорость сходимости для выпуклых функций составляет $O(n^2 \ln (\e^{-1}))$ \cite{nemirovsky1983problem}. Поэтому область применения методов типа эллипсоидов ограничивается размерностью задачи $n \lesssim 100$. С другой стороны, градиентные методы сходятся медленнее относительно точности, но их скорость сходимости не зависит от размерности. К примеру, скорость сходимости быстрого градиентного метода в выпуклом гладком случае составляет $O(\e^{-1/2})$ \cite{nesterov2004introduction}. Таким образом, при $n \gg 100$ градиентные методы выигрывают у методов нулевого порядка.
    
    Наконец, рассмотрим методы высокого порядка. Так как для любой функции построение алгоритма оптимизации, основанного на производных как первого так и более высоких порядков, подразумевает аппроксимацию целевой функции многочленом Тейлора, возникает вопрос о выпуклости полученной аппроксимации. До недавнего времени научное сообщество весьма пессимистично смотрело на данную проблему \cite{baes2009estimate}. Однако, Юрий Нестеров в своей работе \cite{nesterov2021implementable} показал, как можно правильным образом регуляризовать аппроксимацию Тейлора, чтобы сделать её выпуклой. В этой же статье он предложил нижние оценки скорости сходимости тензорных методов для выпуклой оптимизации. Это породило огромный поток работ на эту тему, к примеру \cite{gasnikov2019optimal,bubeck2019near,jiang2019optimal,dvurechensky2019near,ostroukhov2020tensor}. В частности, в работах \cite{gasnikov2019optimal,bubeck2019near,jiang2019optimal} авторы практически одновременно предложили ускоренные тензорные методы со скоростью сходимости $\tilde O \left( \e^{-\frac{2}{3p + 1}} \right)$, где под тильдой в $\tilde O$ обозначен мультипликативный логарифмический фактор. Данная оценка является почти оптимальной и совпадает с нижней оценкой из \cite{nesterov2021implementable}, не считая логарифма. В следующей своей работе \cite{nesterov2021superfast} Нестеров показал, как можно в предположении о Липшицевости производной третьего порядка решать задачу методом второго порядка и предложил "супербыстрый тензорный метод". Полученный таким образом метод имеет сложность, которая оказывается лучше, чем существующие оценки для методов второго порядка. 
    Связяно это с тем, что раньше нижние оценки для алгоритмов $p$-го порядка предлагались, исходя из предположения о Липшицевости также $p$-й производной целевого функционала.
    Таким образом, предположение о Липшицивости производной $(p+1)$-го порядка нас выводит из рассматриваемого класса. Одним из продолжений этой работы послужила статья \cite{ahookhosh2021high_a}. В ней авторы показали, как можно с помощью ускоренных тензорных методов решать задачи с композитом типа 
    \begin{equation*}
        \min_{x \in \text{dom} \psi} \{ F(x) \equiv f(x) + \psi(x)\},
    \end{equation*}
    где $f: \E \to \R$ -- выпуклая замкнутая и, возможно, недифференцируемая функция, $\psi: \E \to \R$ -- простая выпуклая замкнутая функция, $\text{dom} \psi \subseteq \text{int} (\text{dom} \psi)$, $\E$ -- конечномерное вещественное векторное пространство.  Это позволяет, к примеру, решать задачи с ограничениями на простых замкнутых множетствах, если в качестве композита $\psi(x)$ использовать индикатор этого множества.
    
    Очевидно, что скорость сходимости тензорных методов превышает скорость сходимости градиентного метода с точки зрения количества итераций. Но стоимость одной итерации возрастает ввиду необходимости использования производных высокого порядка. Обычно предполагается, что тензорные методы показывают свою эффективность при размерности задачи $n < 1000$. Таким образом, можно сделать вывод, что тензорные методы занимают некую нишу между методами типа эллипсоидов ($n \lesssim 100$) и градиентными методами ($n \gg 100$). В связи с этим, кажется целесообразным исследовать применение методов высокого порядка к задаче \eqref{eq:min-min statement}.
    
    Итак, в исследуемой нами постановке задачи \eqref{eq:min-min statement} мы предполагаем $m \gg n$, $100 < n < 1000$.
    Дополнительно, мы предполагаем, что $Q_x \subset R^m$ -- компактное множество, для $Q_y$ мы рассматриваем 2 случая: $Q_y = \R^n$ и $Q_y \subset \R^n$ -- компакт. $F(x, y)$ является $L_{p, y}$-гладкой по $y$ (см. Определение \ref{fed:L_p-smooth}), $L_{xy}$-гладкой по совокупности переменных (см. Определение \ref{fed:L-smooth by both vars}), $\mu_x$-сильно выпуклой по $x$ и $\mu_y$-сильно выпуклой по $y$. 
    Предлагается, подобно работе \cite{gladin2021solving_b} разбить данную задачу на две подзадачи: внешнюю ($\min_{x \in \R^m}$) и внутреннюю ($\min_{y \in \R^n}$). Внутренняя задача неточно решается супербыстрым тензорным методом. Внешняя задача, используя неточный оракул, полученный из решения внутренней задачи, решается быстрым градиентным методом.
    
    Структура работы выглядит следующим образом. В следующем параграфе приводятся некоторые используемые по ходу статьи обозначения и общие определения. Далее описываются используемые в нашей работе алгоритмы: супербыстрый тензорный метод и быстрый градиентный метод на простых множествах. Наши основные результаты об объединении упомянутых алгоритмов в смешанный оракул приводятся в последующем параграфе. Тут же сравниваются теоретические сложности полученного метода и обычного быстрого градиентного метода. В конце вкратце еще раз описываются полученные нами результаты и обсуждаются возможные дальнейшие направления развития этой темы.
    
\paragraph{Обозначения и определения}\label{par:fo}

    В данном разделе мы опишем некоторые общиеизвестные определения. Также, мы введем используемые обозначения, которые понадобятся в дальнейшем.

    Для некоторой функции $f$ обозначим производную порядка $p$ в точке $x \in \dom f$ по направлениям $h_i \in \R^n, i=1, ..., p$ через $\nabla^p f(x)[h_1,...,h_p],\ p \ge 1$ . Тогда норма $p$-й производной определяется как
    \[
        \|\nabla^p f(x) \|_2 := \max_{h_1, ...,h_p \in \R^n} \{ |\nabla^p f(x)[h_1, ..., h_p]|: \|h_i\| = 1, i = 1, ... p\}
    \]
    или
    \[
        \|\nabla^p f(x) \|_2 := \max_{h \in \R^n} \{|\nabla^p f(x)[h]^p: \|h\|_2 \le 1 \}.
    \]
    
    Обозначим аппроксимацию Тейлора некоторой функции $f$ в точке $\hat x \in \dom f$ вплоть до $p$-го порядка ($p \ge 1$) через
    \begin{gather*}
        f(x) = \Omega_{\hat x, p}^f(x) + o(\|x - \hat x\|^p_2),\ \forall x \in \dom f, \\
        \Omega_{\hat x, p}^f(x) := \sum_{i = 0}^p \frac{1}{i!} \nabla^i f(\hat x)[x - \hat x]^p.
    \end{gather*}
    Также нам понадобится обозначение регуляризованной аппроксимации Тейлора порядка $p \ge 1$:
    \[
        \hat \Omega_{\hat x, p, H}^f(x) := \Omega_{\hat x, p}^f(x) + \frac{H}{(p + 1)!}\|x - \hat x\|_2^{p + 1}, \forall x \in \dom f.
    \]
    Для простоты, если из контекста будет понятно, верхний индекс над $\Omega$ мы будем опускать: $\Omega_{\hat x, p}^f(x) \equiv \Omega_{\hat x, p}(x),\ \hat \Omega_{\hat x, p}^f(x) \equiv \hat \Omega_{\hat x, p}(x)$.
    Заметим, что $\hat \Omega_{\hat x, p, H}^f(x)$ является выпуклой функцией в случае, если $f$ является выпуклой, $L_{p, x}$-гладкой (см. Определение \ref{fed:L_p-smooth}) и $H \ge p L_{p, x}$ \cite{nesterov2021implementable}.

    
    
    
    \begin{fed}\label{fed:L_p-smooth}
        Пусть $f(x)$ -- некоторая $p \ge 1$ раз дифференцируемая функция. Тогда $f$ удовлетворяет условию Липшица $p \ge 1$ порядка (является $L_{p,x}$-гладкой), если
        \begin{equation}\label{eq:L_p-smooth}
            \forall x, x' \in Q_x \Rightarrow \|\nabla^p f(x) - \nabla^p f(x') \|_2 \le L_{p,x} \|x - x'\|_2.
        \end{equation}
    \end{fed}
    В дальнейшем, если будет понятно из контекста, вместо $L_{p,x}$ будем писать $L_p$.
    
    \begin{fed}\label{fed:L-smooth by both vars}
        Пусть $F(x, y)$ -- некоторая дифференцируемая по обеим переменным функция. Тогда $F(x, y)$ удовлетворяет условию Липшица по совокупности переменных (является $L_{xy}$-гладкой по совокупности переменных), если
        \begin{equation}\label{eq:L-smooth by both vars}
            \forall x, x' \in Q_x, y, y' \in Q_y \Rightarrow \|\nabla F(x, y) - \nabla F(x', y') \|_2 \le L_{xy} \|(x, y) - (x', y')\|_2.
        \end{equation}
    \end{fed}
    
    Также, при описании супербыстрого тензорного метода нам понадобится определение дивергенции Брегмана и прокс-функции.
    \begin{fed}
        Дивергенция Брегмана для функции $f(x)$ определяется следующим образом
        \[
            \beta_f(x, y) := f(y) - f(x) - \la \nabla f(x); y - x \ra,\ \forall x, y \in \dom f.
        \]
    \end{fed}
    По сути, дивергенция Брегмана показывает разницу между значением функции в точке $y$ и значением её линейной аппроксимации в точке $y$ относительно точки $x$.
    \begin{fed}
        Прокс-функцией $d_p(x)$ порядка $p$ для некоторого $x \in Q_x$ называется некоторая $p$ раз непрерывно дифференцируемая, 1-сильно выпуклая функция. В нашей работе мы выбираем
        \[
            d_{p}(x) := \frac{1}{p} \|x\|_2^p.
        \]
    \end{fed}
    
    Также, нам понадобятся определения гладкости и сильной выпуклости относительно некоторой функции.
    \begin{fed}
        Пусть $f(x)$ -- некоторая выпуклая функция. Тогда $f(\cdot)$ является $L_h$-гладкой относительно некоторой функции $\rho(\cdot)$, если существует такая константа $L_h > 0$, что $(L_h \rho - h)(\cdot)$ выпукла.
    \end{fed}
    \begin{fed}
        Пусть $f(x)$ -- некоторая выпуклая функция. Тогда $f(\cdot)$ является $\mu_h$-сильно выпуклой относительно некоторой функции $\rho(\cdot)$, если существует такая константа $\mu_h > 0$, что $(h - \mu_h \rho)(\cdot)$ выпукла.
    \end{fed}
    
\paragraph{Используемые алгоритмы}

    В этом разделе мы вкратце описываем используемые методы, опуская многие детали, важные для понимания, и отсылая читателя к указанным первоисточникам.
    
    \subparagraph{Супербыстрый тензорный метод}
    
        Рассмотрим задачу
        \begin{equation}\label{eq:superfast tm problem statement}
            \min_{y \in \R^n} f(y),
        \end{equation}
        где $f$ -- $L_3$-гладкая, $\mu$-сильно выпуклая функция.
    
        Для начала, стоит еще раз акцентировать внимание, на том почему супербыстрый тензорный метод носит такое название. Этот метод использует предположение об $L_3$-гладкости, хотя по факту решает задачу, используя оракул второго порядка. Это стало возможным благодаря замене третьей производной по направлению на её разностную аппроксимацию производными первого порядка (см. Лемма 4.2 в \cite{nesterov2021superfast}). Рассматривая таким образом поставленную задачу, алгоритм получает оценку сходимости лучше, чем нижние оценки для методов второго порядка в традиционной постановке, когда предполагается только $L_2$-гладкость.
        
        В данной работе мы будем использовать ускоренный вариант супербыстрого метода, названного $\text{ATMI}_3$ (Inexact Accelerated 3rd-Order Tensor Method). Перед тем как описывать его псевдокод, опишем подзадачу, которую этот метод решает.
        
        Внутри рассматриваемого метода на каждом шаге необходимо найти неточный минимум функции $\hat \Omega_{y, p, H}(\cdot)$ в следующем смысле. Осуществляется поиск точек во вложенных окрестностях
        \begin{equation}\label{eq:bdgm neighborhood}
            \N_{p, H}^\gamma(y) = \{ T \in \R^n:\ \|\nabla \hat \Omega_{y, p, H}(T)\|_2 \le \gamma \|\nabla f(T)\|_2 \},
        \end{equation}
        где $\gamma \in [0, 1)$ -- некоторый параметр точности. Зафиксируем $\gamma = \frac{1}{2p},\ H = 2p L_p$, и для простоты обозначим
        \[
            \N_p(y) \equiv \N_{p, 2p L_p}^{1 / (2p)}(y).
        \]
        При этом, $p = 3$.
        
        Для решения описанной подзадачи будет использоваться дополнительный алгоритм, который носит название Bregman Distance Gradient Method (BDGM). В \cite{grapiglia2021inexact} было доказано, что BDGM решает данную проблему за линейное время, а в \cite{nesterov2021superfast} было показано, как можно улучшить данный алгоритм, чтобы он работал с неточными градиентами. Данная модификация позволяет избежать вычисления $\nabla^3 f(\hat y_k)[y - \hat y_k]^2,\ \forall \hat y_k, y \in \R^n$, заменив её разностными аналогами с использованием градиентов:
        \begin{equation}
            g_{\hat y_k}^\tau(y) := \frac{1}{\tau^2} (\nabla f(\hat y_k + \tau (y - \hat y_k)) + \nabla f(\hat y_k - \tau(y - \hat y_k)) - 2 \nabla f(\hat y_k)).
        \end{equation}
        Введем дополнительное обозначение $\phi_k(y) \equiv \nabla \hat \Omega_{\hat y_k, 3, 6L_p}(y)$. Тогда аппроксимацию $\phi(y)$ можно переписать в виде
        \begin{equation}\label{eq:approximate subproblem grad}
            g_{\phi_k, \tau}(y) := \nabla f(\hat y_k) + \nabla^2 f(\hat y_k)[y - \hat y_k] + \frac{1}{2} g_{\hat y_k}^\tau(y) + L_3\|y - \hat y_k\|_2^2 (y - \hat y_k).
        \end{equation}
        
        \begin{algorithm}[t]
            \caption{Bregman Distance Gradient Method (BDGM) \cite{nesterov2021superfast}}\label{alg:bdgm}
            \begin{algorithmic}[1]
                \REQUIRE $\delta > 0,\ \hat y_k \in \R^n$.
                
                \STATE $z_0 = \hat y_k$.
                
                \STATE $\tau = \frac{3 \delta}{8(2 + \sqrt 2) \|\nabla f(\hy_k)\|_2}$.
                
                \STATE Определим множество
                \[
                    S_k = \left\{ z: \|z - \hat y_k\| \le 2 \left( \frac{2 + \sqrt 2}{L_3} \|\nabla f(\hat y_k)\|_2 \right)^\frac{1}{3} \right\}.
                \]
                
                \STATE Определим функцию
                \[
                    \rho_k(z) = \frac{1}{2} \left\la \nabla^2 f(\hat y_k)^T (z - \hat y_k); z - \hat y_k \right\ra + L_3 d_4(z - \hat y_k).
                \]
                
                \STATE $k = 0$.
                
                \WHILE{$\|g_{\phi_k, \tau}(z_k)\| > \frac{1}{6} \|\nabla f(z_k)\|_2 - \delta$}
                    
                    \STATE Вычислить $g_{\phi_k, \tau}(z_k)$ через \eqref{eq:approximate subproblem grad}.
                    
                    \STATE 
                    \[
                        z_{k + 1} = \arg \min_{z \in S_k} \left\{ \la g_{\phi_k, \tau}(z_k); z - z_k \ra + 2 \left( 1 + \frac{1}{\sqrt 2} \right) \beta_{\rho_k}(z_k, z) \right\}.
                    \]
                    
                    \STATE $k = k + 1$.
                \ENDWHILE
                
                \RETURN $z_k$.
            \end{algorithmic}
        \end{algorithm}
        
        Псевдокод BDGM можно увидеть в Алгоритме \ref{alg:bdgm}.
        Авторы в своей работе показывают, что при выборе погрешности при решении задачи \eqref{eq:bdgm neighborhood} в виде
        \[
            \delta = O \left( \frac{\e^\frac{3}{2}}{\|\nabla f(\hy_k)\|_2^\frac{1}{2} + \|\nabla^2 f(\hat y_k)\|^\frac{3}{2} / L_3^\frac{1}{2}} \right)
        \]
        Алгоритму \ref{alg:bdgm} нужно будет сделать
        \[
            T_k(\delta) = O \left( \ln \frac{G + H}{\e} \right)
        \]
        итераций, где $G$ и $H$ -- равномерные верхние границы норм градиентов и гессианов, вычисленных в точках, сгенерированных основным алгоритмом. Следовательно, Алгоритму \ref{alg:bdgm} нужно будет один раз вычислить $\nabla^2 f(\cdot)$ и $T_k(\delta)$ раз вычислить $\nabla f(\cdot)$.
        
        \begin{algorithm}[t]
            \caption{Inexact Accelerated 3rd-Order Tensor Method}\label{alg:atmi}
            \begin{algorithmic}[1]
                \REQUIRE $y_0 \in \R^n,\ N \in \mathbb N$
                
                \STATE $c_3 = \left( 5 / (7 L_3) \right)^\frac{1}{3}$
                
                \STATE $\psi_0(y) = d_4(y - y_0)$
                
                \FOR{$i = 0,..., N - 1$}
            
                    \STATE $v_i = \arg \min_{y \in \R^n} \psi_i(y)$.
                
                    \STATE 
                    \[
                        A_i = 2 \left( \frac{2}{3} c_3 \right)^3 \left( \frac{i}{4} \right)^4,\ a_{i + 1} = A_{i + 1} - A_i.
                    \]
                    
                    \STATE
                    \[
                        z_i = \frac{A_i}{A_{i + 1}} y_i + \frac{a_i}{A_{i + 1}} v_i.
                    \]
                    
                    \STATE Вычислить $y_{i + 1} = \N_3(z_i)$ с помощью Алгоритма \ref{alg:bdgm}.
                    
                    \STATE $\psi_{i + 1}(y) = \psi_i(y) + a_{i + 1} \left( f(y_{i + 1}) + \left\la \nabla f(y_{i + 1}); y - y_{i + 1} \right\ra \right)$.
            
                \ENDFOR
                
                \RETURN $y_N$.
            \end{algorithmic}
        \end{algorithm}
        
        Обычный (неускоренный) супербыстрый метод, по сути, представляет из себя Алгоритм \ref{alg:bdgm}, запущенный некоторое предопределенное количество раз. При помощи стандартной техники оценивающих последовательностей этот метод можно ускорить и получить Алгоритм \ref{alg:atmi}.
        Скорость сходимости данного метода приводится в следующей теореме.
        
        \begin{teo}[с. 26 в \cite{nesterov2021superfast}]
            Пусть некоторая функция $f: \R^n \to \R$ является выпуклой и $L_3$-гладкой. Тогда для Алгоритма \ref{alg:atmi} выполняется
            \begin{equation}\label{eq:superfast main convergence rate}
                f(y_N) - f(y^*) \le \frac{7}{60} \left( \frac{6}{N} \right)^4 \cdot L_3 R^4,
            \end{equation}
            где $R = \|y_0 - y^*\|_2$.
            Соответственно, для нахождения $y_\e \in \dom f: f(y_\e) - f(y^*) \le \e$ Алгоритму нужно
            \begin{equation}\label{eq:superfast number of Hessians}
                K = O \left( R \left( \frac{L_3}{\e} \right)^\frac{1}{4} \right)
            \end{equation}
            вычислений гессианов и
            \begin{equation}\label{eq:superfast number of gradients}
                O \left( R \left( \frac{L_3}{\e} \right)^\frac{1}{4} \log \frac{G + H}{\e_g} \right)
            \end{equation}
            вычислений градиентов,
            где $\e_g$ -- нижняя граница для норм всех градиентов, посчитанных во время решения подзадачи.
        \end{teo}
        
        Так как в этой работе мы предполагаем $\mu_y$-сильную выпуклость внутренней задачи, воспользуемся стандартной процедурой рестартов, чтобы получить еще более ускоренный вариант алгоритма.
        
        \begin{algorithm}[t]
            \caption{Restarted Inexact Accelerated 3rd-Order Tensor Method}\label{alg:re-atmi}
            \begin{algorithmic}[1]
                \REQUIRE $\e > 0, y_0 \in \R^n, R \ge \|y_0 - y^*\|_2$
                \FOR{$i = 0, ..., \left\lceil \log \frac{\mu R^2}{\e} \right\rceil - 1$}
                    \STATE $R_i = R/2^i$.
                    \STATE 
                    \[
                        N_i = 6 \left\lceil \sqrt[4]{\frac{7 L_3 R_i^2}{15 \mu}} \right\rceil.
                    \]
                    \STATE Запустить Алгоритм \ref{alg:atmi} в течение $N_i$ итераций, на выходе получить $y_{N_i}$.
                \ENDFOR
            \RETURN $y_{N_i}$
            \end{algorithmic}
        \end{algorithm}
        
        \begin{teo}\label{teo:superfast for strongly convex}
            Пусть решается задача \eqref{eq:superfast tm problem statement}. Тогда для нахождения $y_\e \in \dom f: f(y_\e) - f(y^*) \le \e$ Алгоритму \ref{alg:re-atmi} нужно
            \begin{equation}\label{eq:restarted superfast number of Hessians}
                O \left( \left( \frac{L_3 R^2}{\mu} \right)^\frac{1}{4} \log \frac{\mu R^2}{\e} \right)
            \end{equation}
            вычислений гессианов и 
            \begin{equation}\label{eq:restarted superfast number of gradients}
                O \left( \left( \frac{L_3 R^2}{\mu} \right)^\frac{1}{4} \log \frac{\mu R^2}{\e} \log \frac{G + H}{\e_g} \right)
            \end{equation}
            вычислений градиентов, где $R = \|y_0 - y^*\|_2$, а $\e_g$ -- нижняя граница для норм всех градиентов, посчитанных во время решения подзадачи.
        \end{teo}
        \begin{proof}
            Из \eqref{eq:superfast main convergence rate} и сильной выпуклости получаем
            \[
                \frac{\mu}{2} \|y_N - y^*\|_2^2 \le f(y_N) - f^* \le \frac{7}{60} \left( \frac{6}{N} \right)^4 \cdot L_3 R^4. 
            \]
            Выберем $N_1: \|y_{N_1} - y^*\|_2^2 \le \frac{1}{2} \|y_0 - y^*\|_2^2 = \frac{1}{2} R^2$. Тогда
            \[
                N_1 = 6 \left\lceil \sqrt[4]{\frac{7 L_3 R^2}{15 \mu}} \right\rceil.
            \]
            Далее, аналогично будем выбирать $N_i$, каждый раз уменьшая квадрат расстояния вдвое.
            
            Обозначим $R_k = R / 2^{k - 1}$. Так как мы хотим получить $f(y_{N_k}) - f^* \le \e$, то оценим количество необходимых рестартов:
            \[
                f(y_{N_k}) - f^* \le \frac{7}{60} \left( \frac{6}{N_k} \right)^4 \cdot L_3 R_k^4 \le \frac{\mu R^2}{2^{k - 1}} = \e.
            \]
            В итоге, получаем
            \[
                k = \left\lceil \log \frac{\mu R^2}{\e} \right\rceil + 1.
            \]
            Тогда общее число итераций
            \[
                N = \sum_{i = 1}^k N_i \le kN_1 = O \left( \left( \frac{L_3 R^2}{\mu} \right)^\frac{1}{4} \log \frac{\mu R^2}{\e} \right).
            \]
        \end{proof}\qed
        
    \subparagraph{Неточный проксимальный тензорный метод}
    
        Предположим, что мы решаем задачу с ограничениями на некотором простом замкнутом множестве:
        \begin{equation}
            \min_{y \in Q_y} f(y),
        \end{equation}
        где $f$ -- $L_p$-гладкая, выпуклая функция. Эту задачу также можно переписать в виде
        \begin{equation}\label{eq:composite statement}
            \min_{y \in \dom \psi} \{ \Phi(y) := f(y) + \psi(y) \},
        \end{equation}
        где $\psi(y)$ -- индикаторная функция $Q_y$:
        \[
            \psi(y) = 
            \begin{cases}
                0,& y \in Q_y, \\
                +\infty,& y \not\in Q_y.
            \end{cases}
        \]
        
        
        Введем следующее определение.
        
        \begin{fed}
            Композитным проксимальным оператором $p$-го порядка для некоторой функции $\Phi$ из \eqref{eq:composite statement} называется следующая функция
            \begin{equation}\label{eq:proximal p-th order composite operator}
                \prox_{\Phi/H}^p(\ty) = \arg \min_{y \in E} \left\{ \Phi(y) + \frac{H}{p + 1} \|y - \ty\|_2^{p + 1} \right\}.
            \end{equation}
        \end{fed}
        В \cite{ahookhosh2021high_a} авторы исследуют проксимальные методы высокого порядка через аппроксимацию оператора \eqref{eq:proximal p-th order composite operator} и его неточное решение на каждом шаге в виде подзадачи. Множество возможных решений \eqref{eq:proximal p-th order composite operator} описывается через 
        \begin{equation}\label{eq:proximal p-th order composite operator solutions}
            \mathcal A_H^p(\ty, \gamma) = \left\{ (y, g) \in \dom \psi \times \R^n : g \in \partial \psi(x), \|\nabla f_{\ty, H}^p (y) + g\|_2 \le \gamma \|\nabla f(y) + g\|_2 \right\},
        \end{equation}
        где
        \begin{equation}
            f_{\ty, H}^p (y) = f(y) + H d_{p + 1}(y - \ty),
        \end{equation}
        где $\gamma \in [0, 1)$ -- параметр точности.
        
        Определим функцию
        \begin{equation}\label{eq:rho}
            \rho_{\hy_k, H}(z) = \sum_{k = 1}^{\lfloor p / 2 \rfloor} \frac{1}{(2k)!} D^{2k} f(\hy_k)[z - \hy_k]^{2k} + H d_{p + 1}(z - \hy_k).
        \end{equation}
        Тогда для нахождения решения из \eqref{eq:proximal p-th order composite operator solutions} авторы предлагают Алгоритм \ref{alg:Non-Euclidean Composite Gradient Algorithm}, где для простоты используется обозначение $\rho \equiv \rho_{\hy_k, H}$.
        
        \begin{algorithm}
            \caption{Non-Euclidean Composite Gradient Algorithm \cite{ahookhosh2021high_a}}
            \label{alg:Non-Euclidean Composite Gradient Algorithm}
            \begin{algorithmic}[1]
                \REQUIRE $\hy_k \in \dom \psi$, $\gamma \in [0, 1/p]$, $L > 0$.
                
                \STATE $z_0 = \hy_k$.
                
                \STATE $i = 0$.
                                
                \WHILE{$z_i \not\in \A_H^p(\hy_k, \gamma)$}

                    \STATE Вычислить 
                    \[
                        z_{i + 1} = \arg \min_{z \in \dom \psi} \left\{ \left\la \nabla f_{y_k, H}^p(z_i); z - z_i \right\ra + \psi(z) + 2 L \beta_\rho(z_i, z) \right\}.
                    \]
                    
                    \STATE $g = L(\nabla \rho(z_i) - \nabla \rho(z_{i + 1})) - \nabla f_{y_k, H}^p(z_i),\ g \in \partial \psi(z_{i + 1})$.
                    
                    \STATE $i = i + 1$.
                
                \ENDWHILE
                
                \RETURN $(z_i, g)$.
            \end{algorithmic}
        \end{algorithm}
        
        Согласно Теореме 3.4 из \cite{ahookhosh2021high_a}, данный алгоритм сходится линейно относительно точности решения основной задачи.
        
        Далее, для решения исходной задачи \eqref{eq:composite statement} авторы оборачивают Алгоритм \ref{alg:Non-Euclidean Composite Gradient Algorithm} в итерационный процесс и ускоряют с помощью техники оценивающих последовательностей. Псевдокод полученного метода приведен в Алгоритме \ref{alg:Bi-Level High-Order Algorithm}. Его сложность описана в следующей теореме.
        
        \begin{algorithm}
            \caption{Bi-Level High-Order Algorithm}
            \label{alg:Bi-Level High-Order Algorithm}
            \begin{algorithmic}[1]
                \REQUIRE $y_0 \in \dom \psi,\ \gamma \in [0, 1/p]$
                
                \STATE $H = \frac{6}{(p - 1)!} L_p,\ A_0 = 0,\ c_p = \left( \frac{1 - \gamma}{H} \right)^\frac{1}{p}$.
                
                \STATE $\Psi_0(y) = d_{p + 1}(y - y_0)$
                
                \STATE $k = 0$
                
                \WHILE{$\Phi(y_k) - \Phi(y^*) > \e$}
                
                    \STATE $v_k = \arg \min_{x \in \dom \psi} \Psi_k(x)$
                    
                    \STATE
                    \[
                        A_k = \left( \frac{c_p}{2} \right)^p \left( \frac{k}{p + 1} \right)^{p + 1},\ a_{k + 1} = A_{k + 1} - A_k.
                    \]
                    
                    \STATE
                    \[
                        z_k = \frac{A_k}{A_{k + 1}} y_k + \frac{a_{k + 1}}{A_{k + 1}} v_k.
                    \]
                    
                    \STATE Вычислить $(T_k, g) \in \A_H^p(z_k, \gamma)$ с помощью Алгоритма \ref{alg:Non-Euclidean Composite Gradient Algorithm}.
                    
                    \STATE Найти $y_{k + 1}: \Phi(y_{k + 1}) \le \Phi(T_k)$.
                    
                    \STATE $\Psi_{k + 1}(y) = \Psi_k(x) + a_{k + 1}(\Phi(y_{k + 1}) + \la \nabla \Phi(y_{k + 1}); y - y_{k + 1} \ra)$.
                    
                    \STATE $k = k + 1$.
                
                \ENDWHILE
                
                \RETURN $y_k$.
            \end{algorithmic}
        \end{algorithm}
        
        \begin{teo}[Теорема 3.11 в \cite{ahookhosh2021high_a}]
            Пусть решается задача \eqref{eq:composite statement}, $H \ge p L_p$, $p \ge 2,\ q = \lfloor p / 2 \rfloor$, $\gamma \in [0, 1/p]$, $\rho$ определяется из \eqref{eq:rho}. Тогда сложность Алгоритма \ref{alg:Bi-Level High-Order Algorithm} определяется как
            \begin{equation}
                \tilde O \left( \frac{L_p R^{p + 1}}{\e} \right)^\frac{1}{p + 1},
            \end{equation}
            где $R = \|y_0 - y^*\|_2$.
        \end{teo}
        
        Дополнительно предположив сильную выпуклость $f(x)$, мы можем улучшить полученный результат с помощью техники рестартов аналогично Теореме \ref{teo:superfast for strongly convex}. Как итог, мы получим следующий результат.
        
        \begin{teo}\label{teo:proximal tensor method for strongly convex}
            Пусть решается задача \eqref{eq:composite statement}, $H \ge p L_p$, $p \ge 2,\ q = \lfloor p / 2 \rfloor$, $\gamma \in [0, 1/p]$, $\rho$ определяется из \eqref{eq:rho}. Предположим также, что $f(x)$ является $\mu$-сильно выпуклой. Тогда, применив технику рестартов к Алгоритму \ref{alg:Bi-Level High-Order Algorithm}, получим метод со сложностью
            \begin{equation}
                O \left( \left( \frac{L_p R^{p - 1}}{\mu} \right)^\frac{1}{p + 1} \log \frac{\mu R^2}{\e} \right)
            \end{equation}
        \end{teo}
        \begin{proof}
            Доказывается аналогично Теореме \ref{teo:superfast for strongly convex}.
        \end{proof}\qed
        
    \subparagraph{Быстрый градиентный метод на простом множестве}
    
        Рассмотрим следующую постановку задачи
        \begin{equation}\label{eq:fgm on simple set problem statement}
            \min_{x \in Q_x} f(x),
        \end{equation}
        где $Q_x \subset \R^m$ -- некоторое простое множество, $f$ -- $L$-гладкая $\mu$-сильно выпуклая функция. Под простым множеством подразумевается такое множество, проекцию градиента на которое можно явно вычислить. Определим градиентное отображение на $Q_x$:
        \begin{fed}[Определение 2.2.3 в \cite{nesterov2004introduction}]
            Зафиксируем некоторое $\gamma > 0$. Обозначим
            \begin{align*}
                x_{Q_x}(\hat x, \gamma) &:= \arg \min_{Q_x} f(\hat x) + \la \nabla f(\hat x); x - \hat x \ra + \frac{\gamma}{2} \|x - \hat x\|_2^2, \\
                g_{Q_x}(\hat x, \gamma) &:= \gamma (\hat x - x_{Q_x}(\hat x, \gamma).
            \end{align*}
            Тогда $g_{Q_x}(\hat x, \gamma)$ называется градиентным отображением $f$ в точке $\hat x$ на множество $Q_x$.
        \end{fed}
        
        Псевдокод быстрого градиентного метода для задачи \eqref{eq:fgm on simple set problem statement} приведен в Алгоритме \ref{alg:fgm on simple set}. Сложность данного метода приведена в следующей теореме.
        
        \begin{algorithm}[t]
            \caption{Быстрый градиентный метод на простом множестве}\label{alg:fgm on simple set}
            \begin{algorithmic}[1]
                \REQUIRE $x_0 \in Q_x$, размер шага $\alpha \in (0, 1)$, $N \in \mathbb N$.
                    
                \STATE $y_0 = x_0,\ q = \mu / L$
                    
                \FOR{$i = 1, ..., N - 1$}
                
                    \STATE $x_{i + 1} = x_{Q_x}(y_k, L)$
                    
                    \STATE Вычислить $\alpha_{i + 1} \in (0, 1): \alpha_{i + 1}^2 = (1 - \alpha_{i + 1})\alpha_i^2 + q \alpha_{i + 1}.$
                    
                    \STATE $\beta_i = \frac{\alpha_i(1 - \alpha_i)}{\alpha_i^2 + \alpha_{i + 1}}$
                    
                    \STATE $y_{i + 1} = x_{i + 1} + \beta_k(x_{i + 1} - x_i)$.
                    
                \ENDFOR
                
                \RETURN $x_{N}$
            \end{algorithmic}
        \end{algorithm}
        
        \begin{teo}[с. 119 в \cite{nesterov2004introduction}]
            Пусть решается задача \eqref{eq:fgm on simple set problem statement}. Тогда для нахождения $x_\e \in Q_x: f(x_\e) - f(x^*)$ Алгоритму \ref{alg:fgm on simple set} необходимо
            \begin{equation}\label{eq:fgm on simple set convergence rate}
                O \left( \frac{L R^2}{\e} \right)^\frac{1}{2}
            \end{equation}
            итераций в выпуклом случае и
            \begin{equation}\label{eq:fgm on simple set convergence rate strongly convex}
                O \left( \left( \frac{L}{\mu} \right)^\frac{1}{2} \log \frac{L R^2}{\e} \right)
            \end{equation}
            в сильно выпуклом случае, где $R = \|x_0 - x^*\|_2$.
        \end{teo}
        
\paragraph{Полученные результаты}

    Опишем еще раз рассматриваемую нами постановку задачи.
    \begin{equation}\label{eq:min-min statement last}
        \min_{x \in Q_x} \min_{y \in Q_y} F(x, y),
    \end{equation}
    где выполняются следующие предположения
    \begin{enumerate}
        \item $Q_x \subset \R^m$ -- выпуклое компактное множество, $m \gg 100$;
        
        \item $Q_y = \R^n$ или $Q_y \subset \R^n$ -- компакт, $100 < n < 1000$;
        
        \item $F(x, y)$ -- $\mu_y$-сильно выпуклая по $y$;
        
        \item $F(x, y)$ -- $\mu_x$-сильно выпуклая по $x$;
        
        \item $F(x, y)$ -- $L_{p, y}$-гладкая по $y$;
        
        
        \item $F(x, y)$ -- $L_{xy}$-гладкая по совокупности переменных.
    \end{enumerate}
    
    Как уже упоминалось ранее, для решения задачи \eqref{eq:min-min statement last} мы используем смешанный оракул. 
    На первом этапе мы неточно решаем внутреннюю задачу: для некоторого фиксированного $x_{k} \in Q_x$ находим
    \begin{equation}\label{eq:inner problem}
        \tilde y_{k + 1} = \ty(x_k) \in Q_y: F(x_k, \ty_{k + 1}) - F(x_k, y^*_{k + 1}) \le \e,
    \end{equation}
    где $y^*_{k + 1} = y^*(x_{k + 1}) = \min_{y \in Q_y} F(x_k, y)$. Упомянутые предположения о функции позволяют нам использовать для решения супербыстрые методы для сильно выпуклых функций (см. Теорему \ref{teo:superfast for strongly convex}), если мы решаем внутреннюю задачу на всём пространстве. Если же мы решаем внутреннюю задачу на компакте, то мы можем использовать проксимальный тензорный метод для задачи с композитом (см. Теорему \ref{teo:proximal tensor method for strongly convex}).
    На следующем этапе, исходя из полученного решения внутренней задачи, мы осуществляем итерацию внешнего метода для решения внешней задачи и на выходе получаем $x_{k + 1}$. Ввиду предположений о \eqref{eq:min-min statement last}, здесь мы используем быстрый градиентный метод для сильно выпуклых функций на простых множествах. Далее мы вновь решаем внутреннюю задачу уже при фиксированном $x_{k + 1}$ и т. д.
    
    Так как мы решаем внутреннюю задачу неточно, необходимо эту неточность учитывать при обращении к оракулу во время решения внешней задачи. Для этого мы используем концепцию $(\delta, L)$-оракула.
    
    \begin{fed}[Определение 1 в \cite{gasnikov2015universal}]
        $(\delta, L)$-оракул выдает на запрос, в котором указана одна точка $x \in Q_x$ такие ($f_\delta(x), g_\delta(x)$), что 
        \[
            \forall x' \in Q_x \Rightarrow 0 \le f(x') - f_\delta(x) - \la g_\delta(x), x - x' \ra \le \frac{L}{2} \|x' - x\|^2_2 + \delta.
        \]
    \end{fed}
    
    Прежде чем описать, как можно получить неточный оракул из неточного решения внутренней задачи, введем вспомогательное Предложение.
    
    \begin{pro}[Утверждение 3 в \cite{gasnikov2015universal}]\label{pro:from gasnikov}
        Пусть $F: Q_x \times Q_y \to \R$, $f(x) = \min_{y \in Q_y} F(x, y)$, где $F(x, y)$ -- выпуклая, $L_{xy}$-гладкая по совокупности переменных функция, а $Q_x \subseteq \R^m$, $Q_y \subseteq \R^n$ -- некоторые выпуклые множества. Пусть 
        \begin{equation}\label{eq:from gasnikov, main assumption}
            \forall x \in Q_x\ \exists y_\delta = y_\delta(x) \in Q_y: \max_{y' \in Q_y} \la F(x, y_\delta); y_\delta - y' \ra \le \delta.
        \end{equation}
        Тогда
        \begin{gather*}
            \forall x' \in Q_x \Rightarrow \|\nabla f(x') - \nabla f(x)\|_2 \le L_{xy} \|x' - x\|_2, \\
        \end{gather*}
        и $(F(x, y_\delta) - 2\delta, \nabla_x F(x, y_\delta))$ будет $(6 \delta, 2L_{xy})$-оракулом для $f(x)$ на $Q_x$.
    \end{pro}
    
    Как видно, для выполнения этого Предложения нам необходимо выполнение условия \eqref{eq:from gasnikov, main assumption}. Так как мы рассматриваем внутреннюю задачу на двух разных доменах, оба случая нам нужно рассмотреть по отдельности.
    
    \subparagraph{Внутренняя задача на всём пространстве}

        \begin{pro}\label{pro:anti-linear approx upper bound}
            Пусть $f: \R^n \to \R$ -- некоторая $\mu$-сильно выпуклая $L$-гладкая функция. Пусть $\exists x_\e \in \R^n: f(x_\e) - f(x^*) \le \e$, где $x^* = \arg \min_{x \in \R^n} f(x)$. Тогда 
            \[
                \forall x \in \R^n \Rightarrow \la \nabla f(x_\e); x_\e - x \ra \le L \|x_\e - x\|_2 \sqrt \frac{2 \e}{\mu}.
            \]
        \end{pro}
        \begin{proof}
            По неравенству Коши-Буняковского-Шварца
            \begin{equation}\label{eq:pro_1 cauchi-bun-schw}
                \la \nabla f(x_\e); x_\e - x \ra \le \|\nabla f(x_\e)\|_2 \|x_\e - x\|_2.
            \end{equation}
            Из условия Липшица \eqref{eq:L-smooth by both vars}
            \[
                \|\nabla f(x_\e)\|_2 = \|\nabla f(x_\e) - \nabla f(x^*)\|_2 \le L \|x_\e - x^*\|_2.
            \]
            Далее, из сильной выпуклости и определения $x_\e$ имеем
            \[
                \frac{\mu}{2} \|x_\e - x^*\|^2_2 \le f(x_\e) - f(x^*) \le \e \Leftrightarrow \|x_\e - x^*\|_2 \le \sqrt \frac{2 \e}{\mu}.
            \]
            Таким образом,
            \[
                \|\nabla f(x_\e)\|_2 \le L \sqrt \frac{2 \e}{\mu}.
            \]
            Подставляя данное неравенство в \eqref{eq:pro_1 cauchi-bun-schw}, получаем искомый результат.
        \end{proof}\qed
        
        Очевидно, что для выполнения Предложения \ref{pro:from gasnikov} в условиях нашей задачи необходимо оценить $\|x_\e - x\|_2$ в Предложении \ref{pro:anti-linear approx upper bound}. Мы сделаем это в следующей теореме, аналогичной Теореме 4 в \cite{gladin2021solving_a}.
    
        \begin{teo}\label{teo:about (delta, L) oracle and smoothness}
            Пусть $F: Q_x \times \R^n \to \R$, $f(x) = \min_{y \in \R^n} F(x, y)$. При этом, 
            \begin{itemize}
                \item $Q_x \subset \R^m$ -- некоторое компактное множество,
                
                \item $F(x, y)$ выпуклая по $x$,
                
                \item $F(x, y)$ $\mu_y$-сильно выпуклая по $y$,
                
                \item $F(x, y)$ $L_y$-гладкая по $y$,
                
                \item $F(x, y)$ $L_{xy}$-гладкая по совокупности переменных функция.
            \end{itemize}
            Тогда $\forall x \in Q_x$ $\exists y_\e = y_\e(x) \in \R^n:$
            \begin{enumerate}
                \item 
                \[
                    F(x, y_\e) - f(x) \le \e,
                \]
                
                \item $\forall x' \in Q_x \Rightarrow$
                \[
                    \|\nabla f(x') - \nabla f(x)\|_2 \le L_{xy} \|x' - x\|_2
                \]
                
                \item $(F(x, y_\e) - 2\delta, \nabla_x F(x, y_\e))$ будет $(6 \delta, 2L_{xy})$-оракулом для $f(x)$ на $Q_x$, где 
                \[
                    \delta = \frac{2L_y}{\mu_y} \left(\e + 2 \sqrt{D \e} \right),\ D = \max_{x \in Q_x} \left( F(x, 0) - F(x, y(x)) \right).
                \]
            \end{enumerate}
        \end{teo}
        \begin{proof}
            Из Предложения \ref{pro:anti-linear approx upper bound} мы знаем, что
            \begin{gather}
                \forall x \in Q_x,\ \forall y = y(x) \in \R^n\ \exists y_\e = y_\e(x) \in \R^n: \notag \\ 
                \la \nabla_y F(x, y_\e); y_\e - y \ra \le L_y \|y_\e - y\|_2 \sqrt{\frac{2 \e}{\mu_y}}. \label{eq:from pro_1}
            \end{gather}
            Из неравенства треугольника
            \[
                \|y_\e - y\|_2 \le \|y_\e\|_2 + \|y\|_2 \le \|y_\e\|_2 + \max_{x \in Q_x} \|y(x)\|_2.
            \]
            
            Сначала оценим $\max_{x \in Q_x} \|y(x)\|_2$. Из сильной выпуклости по $y$ и снова из неравенства треугольника
            \[
                \|y(x)\|_2 - \|y'\|_2 \le \|y(x) - y'\|_2 \le \sqrt{\frac{2}{\mu_y} \left( F(x, y') - F(x, y(x)) \right)}.
            \]
            Так как мы можем взять любой $y' \in \R^n$, предположим, что $y' = 0$. Тогда
            \[
                \|y(x)\|_2 \le \sqrt{\frac{2}{\mu_y} \left( F(x, 0) - F(x, y(x)) \right)}.
            \]
            Обозначим $\Delta_0(x) = F(x, 0) - F(x, y(x))$. Так как $\Delta_0(x)$ непрерывна на компакте $Q_x$, то она ограничена. Обозначим $D = \max_{x \in Q_x} \Delta_0(x)$. Тогда мы получаем
            \begin{equation}\label{eq:max of y upper bound}
                \max_{x \in Q_x} \|y(x)\|_2 \le \sqrt{\frac{2}{\mu_y} D}.
            \end{equation}
            
            Теперь оценим $\|y_\e\|_2$. Опять из сильной выпуклости и неравенства треугольника имеем
            \[
                \|y_\e\|_2 - \|y^*\|_2 \le \|y_\e - y^*\|_2 \le \sqrt{\frac{2}{\mu_y} \left( F(x, y_\e) - f(x) \right)} \le \sqrt{\frac{2}{\mu_y} \e}.
            \]
            Так как 
            \[
                \|y^*\|_2 = \|y^*(x)\|_2 \le \max_{x \in Q_x} \|y(x)\|_2 \le \sqrt{\frac{2}{\mu_y} D},
            \]
            то
            \begin{equation}\label{eq:y_eps upper bound}
                \|y_\e\|_2 - \sqrt{\frac{2}{\mu_y} D} \le \|y_\e\|_2 - \|y^*\|_2 \le \sqrt{\frac{2}{\mu_y} \e} \Leftrightarrow \|y_\e\|_2 \le \sqrt{\frac{2}{\mu_y} \e} + \sqrt{\frac{2}{\mu_y} D}.
            \end{equation}
            
            В конце концов, из \eqref{eq:from pro_1}, \eqref{eq:max of y upper bound}, \eqref{eq:y_eps upper bound} получаем конечную оценку
            \[
                \la \nabla_y F(x, y_\e); y_\e - y \ra \le \left( \sqrt{\frac{2}{\mu_y} \e} + 2\sqrt{\frac{2}{\mu_y} D} \right) L_y \sqrt\frac{2 \e}{\mu_y} = \frac{2L_y}{\mu_y} \left( \e + 2 \sqrt{D\e} \right).
            \]
            
            Осталось применить Предложение \ref{pro:from gasnikov}, и мы получаем искомый результат.
        \end{proof}\qed
        
    \subparagraph{Внутренняя задача на компакте}
    
        \begin{pro}\label{pro:anti-linear approx upper bound for compact}
            Пусть решается задача \eqref{eq:composite statement} с метода из Теоремы \ref{teo:proximal tensor method for strongly convex}. При этом, $f(y)$ -- $\mu$-сильно выпуклая, $Q_y \subset \R^n$ -- некоторый компакт. Пусть $\exists y_\e \in \R^n: f(y_\e) - f(y^*) \le \e$, где $y^* = \arg \min_{y \in Q_y} f(y)$. Тогда
            \begin{equation}
                \forall y \in Q_y \Rightarrow \la \nabla \Phi(y_\e); y_\e - y \ra \le \frac{H}{1 - \gamma} \left( \sqrt{\frac{2}{\mu} \e} + \sqrt{\frac{2}{\mu} D} \right)^{p + 1},
            \end{equation}
            где $D := \max_{y \in Q_y} \Phi(y) - \Phi(y^*)$.
        \end{pro}
        \begin{proof}
            Из неравенства Коши-Буняковского-Шварца
            \begin{equation}\label{eq:pro_2 cauchi-bun-schw}
                \la \nabla \Phi(y_\e); y_\e - y \ra \le \|\nabla \Phi(y_\e)\|_2 \|y_\e - y\|_2.
            \end{equation}
            
            В первую очередь, оценим $\|\nabla \Phi(y_\e)\|_2$. По определению,
            \[
                \|\nabla \Phi(y_\e)\|_2 = \|\nabla f(y_\e) + g\|_2,\ g \in \partial \psi(y_\e).
            \]
            Так как $\Phi(y_\e) \le \Phi(T_{k - 1})$, $(T_{k - 1}, g) \in \A_H^p(y_{k - 1}, \gamma)$ в
            Алгоритме \ref{alg:Bi-Level High-Order Algorithm}, то по Лемме 2.2 в \cite{ahookhosh2021high_a} имеем
            \begin{equation}\label{eq:grad phi upper bound}
                \|\nabla \Phi(y_\e)\|_2 = \|\nabla f(y_\e) + g\|_2 \le \frac{H}{1 - \gamma} \|y_\e - y_{k - 1}\|_2^p.
            \end{equation}
            Оценим $\|y_\e - y_k\|_2$. Из неравенства треугольника
            \[
                \|y_\e - y_{k - 1}\|_2 \le \|y_\e - y^*\|_2 + \|y_{k - 1} - y^*\|_2.
            \]
            Так как $f(y)$ $\mu$-сильно выпукла, $\psi(y)$ выпукла, то $\Phi(y)$ $\mu$-сильно выпукла. и для нее выполняется
            \begin{equation*}
                \|y_\e - y^*\|_2 \le \sqrt{\frac{2}{\mu} \left( \Phi(y_\e) - \Phi(y^*) - \la \nabla \Phi(y^*); y_\e - y^* \ra \right)}.
            \end{equation*}
            Из определения $y_\e$ и 
            \begin{equation}\label{eq:linear approximation in minimum for Q_y ge 0}
                \forall y \in Q_y \Rightarrow \la \nabla \Phi(y^*); y - y^* \ra \ge 0
            \end{equation}
            получаем оценку
            \begin{equation}\label{eq:anti-linear approx upper bound for compact, a}
                \|y_\e - y^*\|_2 \le \sqrt{\frac{2}{\mu} \e}.
            \end{equation}
            Далее, снова из сильной выпуклости получаем
            \[
                \|y_{k - 1} - y^*\|_2 \le \sqrt{\frac{2}{\mu} \left( \Phi(y_{k - 1}) - \Phi(y^*) - \la \nabla \Phi(y^*); y_{k - 1} - y^* \ra \right)}.
            \]
            Обозначим $\Delta(y) = \Phi(y) - \Phi(y^*)$. Так как $\Delta(y)$ непрерывна на компакте $Q_y$, то она ограничена на нем, то есть $\exists D \in \R: D = \max_{y \in Q_y} \Delta(y)$. Тогда
            \[
                \forall y \in Q_y \Rightarrow \Phi(y) - \Phi(y^*) \le D.
            \]
            Из этого факта и \eqref{eq:linear approximation in minimum for Q_y ge 0} получаем
            \begin{equation}\label{eq:anti-linear approx upper bound for compact, b}
                \|y_{k - 1} - y^*\|_2 \le \sqrt{\frac{2}{\mu} D}.
            \end{equation}
            Таким образом, из \eqref{eq:grad phi upper bound}, \eqref{eq:anti-linear approx upper bound for compact, a}, \eqref{eq:anti-linear approx upper bound for compact, b} получаем
            \begin{equation}\label{eq:grad phi upper bound new}
                \|\nabla \Phi(y_\e)\|_2 \le \frac{H}{1 - \gamma} \left( \sqrt{\frac{2}{\mu} D} + \sqrt{\frac{2}{\mu} \e} \right)^p.
            \end{equation}
            
            Чтобы оценить $\|y_\e - y\|_2$, также воспользуемся неравенством треугольника, \eqref{eq:linear approximation in minimum for Q_y ge 0}, \eqref{eq:anti-linear approx upper bound for compact, a}, \eqref{eq:anti-linear approx upper bound for compact, b} и получим
            \begin{equation}\label{eq:pro_2 distance upper bound}
                \|y_\e - y\|_2 \le \sqrt{\frac{2}{\mu} D} + \sqrt{\frac{2}{\mu} \e}.
            \end{equation}
            
            Итого, из \eqref{eq:pro_2 cauchi-bun-schw}, \eqref{eq:grad phi upper bound new}, \eqref{eq:pro_2 distance upper bound} получаем
            \[
                \la \nabla \Phi(y_\e); y_\e - y \ra \le \frac{H}{1 - \gamma} \left( \sqrt{\frac{2}{\mu} D} + \sqrt{\frac{2}{\mu} \e} \right)^{p + 1}.
            \]
        \end{proof}\qed
        
        Как видно, в данном случае у нас сразу выполняется условие, необходимое нам в 
        Предложении \ref{pro:from gasnikov}. Так что мы можем привести следующую теорему без доказательства.
        
        \begin{teo}\label{teo:about (delta, L) oracle and smoothness on compact}
            Пусть $F: Q_x \times Q_y \to \R$, $f(x) = \min_{y \in \R^n} F(x, y)$. При этом, 
            \begin{itemize}
                \item $Q_x \subset \R^m$ и $Q_y \subset \R^n$ -- некоторые компактные множества,
                
                \item $F(x, y)$ выпуклая по $x$,
                
                \item $F(x, y)$ $\mu_y$-сильно выпуклая по $y$,
                
                \item $F(x, y)$ $L_y$-гладкая по $y$,
                
                \item $F(x, y)$ $L_{p, y}$-гладкая по $y$,
                
                \item $F(x, y)$ $L_{xy}$-гладкая по совокупности переменных.
            \end{itemize}
            Если решать внутреннюю задачу с помощью метода из Теоремы \ref{teo:proximal tensor method for strongly convex}, то $\forall x \in Q_x$ $\exists y_\e = y_\e(x) \in Q_y:$
            \begin{enumerate}
                \item 
                \[
                    F(x, y_\e) - f(x) \le \e,
                \]
                
                \item $\forall x' \in Q_x \Rightarrow$
                \[
                    \|\nabla f(x') - \nabla f(x)\|_2 \le L_{xy} \|x' - x\|_2
                \]
                
                \item $(F(x, y_\e) - 2\delta, \nabla_x F(x, y_\e))$ будет $(6 \delta, 2L_{xy})$-оракулом для $f(x)$ на $Q_x$, где 
                \[
                    \delta = \frac{H}{1 - \gamma} \left( \sqrt{\frac{2}{\mu_y} D} + \sqrt{\frac{2}{\mu_y} \e} \right)^{p + 1},\ D = \max_{y \in Q_y} \left( F(x, y(x)) - F(x, y^*(x)) \right),
                \]
                где $H \ge p L_{p, y}$, $\gamma \in [0, 1/p]$.
            \end{enumerate}
        \end{teo}
        
    \subparagraph{Внешняя задача}
    
        Итак, мы получили, что решение внутренней задачи с точностью $\e$ дает нам возможность использовать $(6 \delta, 2 L_{xy})$-оракул для внешней задачи.
        В связи с этим можно для внешней задачи применять быстрый градиентный метод. Скорость сходимости быстрого градиентного метода для выпуклой функции в условиях $(\delta, L)$-оракула описывается следующей теоремой.
        
        \begin{teo}[Теорема 4.9 в \cite{devolder2013exactness}]
            Пусть функция $f$ выпукла и наделена $(\delta, L)$-оракулом. Тогда последовательность $x_k$, которая генерируется быстрым градиентным методом с использованием этого оракула, удовлетворяет
            \begin{equation}\label{eq:fgm with inexact oracle complexity}
                f(x_k) - f^* \le \frac{2 L R^2}{(k + 1)^2} + \frac{1}{3}(k + 3) \delta,
            \end{equation}
            где $R = \|x_0 - x^*\|_2$.
        \end{teo}
        
        Псевдокод быстрого градиентного метода с неточным оракулом совпадает с Алгоритмом \ref{alg:fgm on simple set}, где надо заменить выходы оракула первого порядка на их неточные аналоги. Чтобы учесть $\mu_x$-сильную выпуклость $F(x, y)$ можно воспользоваться методом рестартов.
        
        \begin{teo}\label{teo:fgm with inexact oracle for strongly convex}
            Пусть функция $f: Q_x \to \R$, где $Q_x \subset \R^m$ -- выпуклое компактное множество, $\mu$-сильно выпукла и наделена $(\delta, L)$-оракулом. Предположим, что мы можем варьировать $\delta$. Тогда из Алгоритма \ref{alg:fgm on simple set} с неточным оракулом при $q = 0$ с помощью рестартов можно получить алгоритм со сложностью
            \begin{equation}\label{eq:fgm with inexact oracle for strongly convex complexity}
                O \left( \sqrt \frac{L}{\mu} \log \left( \frac{\mu R^2}{\e - \left( \sqrt \frac{10 L}{\mu} + 3 \right) \delta} \right) \right),
            \end{equation}
            где $R = \|x_0 - x^*\|_2$, а $\delta $ выбирается из условия
            \[
                \delta \le \frac{\mu R^2}{10 \left( \sqrt \frac{10 L}{\mu} + 3 \right)}.
            \]
        \end{teo}
        \begin{proof}
            Использовать технику рестартов в её традиционном виде нам мешает аддитивный член в \eqref{eq:fgm on simple set convergence rate}, линейно зависящий от $\delta$. В первую очередь, избавимся от него.
            Так как мы имеем возможность выбирать $\delta$, выберем его таким что
            \begin{equation}\label{eq:delta estimation}
                \frac{1}{\mu} (N + 3) \delta \le \frac{L R^2}{\mu N^2} \Leftrightarrow \delta \le \frac{L R^2}{N^2 (N + 3)}.
            \end{equation}
            Тогда мы можем получить следующую оценку из \eqref{eq:fgm with inexact oracle complexity}
            \[
                \frac{4 L R^2}{\mu N^2} + \frac{1}{\mu} (N + 3) \delta \le \frac{5 L R^2}{\mu N^2}.
            \]
            
            Далее, применяем традиционную технику рестартов, аналогично доказательству Теоремы \ref{teo:superfast for strongly convex}. В итоге, получим $N_i = \sqrt{10L / \mu}, \forall i = 1, ..., k$
            \[
                N = kN_1 =  O \left( \sqrt \frac{L}{\mu} \log \left( \frac{\mu R^2}{\e - \left( \sqrt \frac{10 L}{\mu} + 3 \right) \delta} \right) \right),
            \]
            где $\delta$ выбирается, исходя из \eqref{eq:delta estimation}.
            
        \end{proof}\qed
        
        Можно заметить, что в \eqref{eq:fgm with inexact oracle complexity} присутствует линейное накопление ошибки, возникающей из-за неточности $\delta$. Однако, в \eqref{eq:fgm with inexact oracle for strongly convex complexity}, помимо лучшей скорости сходимости, эта проблема отсутствует. Неточность содержится под логарифмом, в связи с чем пренебрежимо мала.
        
        Таким образом, если мы применим для решения задачи \eqref{eq:min-min statement last} смешанный оракул, который описан в начале параграфа, то получим следующий результат для внутренней задачи всём пространстве.
        
        \begin{teo}\label{teo:main theorem on R^n}
            Пусть $F: Q_x \times \R^n \to \R$, $f(x) = \min_{y \in \R^n} F(x, y)$. Обозначим за $\tilde \e$ точность решения внутренней задачи \eqref{eq:inner problem}: $F(x, y_{\tilde \e}) - f(x) \le \tilde \e$. При этом, 
            \begin{itemize}
                \item $Q_x \subset \R^m$ -- некоторое компактное множество,
                
                \item $F(x, y)$ $\mu_x$-сильно выпуклая по $x$,
                
                \item $F(x, y)$ $\mu_y$-сильно выпуклая по $y$,
                
                \item $F(x, y)$ $L_y$-гладкая по $y$,
                
                \item $F(x, y)$ $L_{3,y}$-гладкая по $y$,
                
                \item $F(x, y)$ $L_{xy}$-гладкая по совокупности переменных функция,
                
                \item обозначим $\delta(\tilde \e) = \frac{2 L_y}{\mu_y} \left( \tilde \e + 2 \sqrt{D \tilde \e} \right),\ D = \max_{x \in Q_x} (F(x, 0) - F(x, y(x)))$ и выберем $\tilde \e$ таким образом, что 
                \[
                    \delta(\tilde \e) \le \frac{\mu_x \|x_0 - x^*\|_2^2}{10 \left( \sqrt \frac{20 L_{xy}}{\mu_x} + 3 \right)}.
                \]
            \end{itemize}
            Будем решать внутреннюю задачу \eqref{eq:inner problem} супербыстрым методом для сильно выпуклых функций (Теорема \ref{teo:superfast for strongly convex}). Внешнюю задачу будем решать быстрым градиентным методом с неточным оракулом для сильно выпуклых функций (Теорема \ref{teo:fgm with inexact oracle for strongly convex}). Тогда сложность конечного алгоритма составит
            \begin{equation}\label{eq:final algorithm complexity}
                O \left( 
                    \left( \frac{L_{xy}}{\mu_x} \right)^\frac{1}{2} \left(  
                        \frac{L_{3,y} R_y^2}{\mu_y}
                    \right)^\frac{1}{4}
                    \log \frac{\mu_y R_y^2}{\tilde \e} \log \frac{G + H}{\e_g} \log \left( \frac{\mu R_x^2}{\e - \left( \sqrt \frac{10 L}{\mu} + 3 \right) \delta} \right)
                \right)
            \end{equation}
            или
            \begin{equation}\label{eq:final algorithm complexity simple}
                \tilde O \left( 
                    \left( \frac{L_{xy}}{\mu_x} \right)^\frac{1}{2} \left(  
                        \frac{L_{3,y} R_y^2}{\mu_y}
                    \right)^\frac{1}{4}
                \right)
            \end{equation}
            где $R_x = \|x_0 - x^*\|_2,\ R_y = \|y_0 - y^*\|_2$, а под тильду в $\tilde O$ мы занесли все логарифмические факторы.
        \end{teo}
        \begin{proof}
            Следует непосредственно из Теорем \ref{teo:superfast for strongly convex}, \ref{teo:about (delta, L) oracle and smoothness}, \ref{teo:fgm with inexact oracle for strongly convex}.
        \end{proof}\qed
        
        И для внутренней задачи, определенной на компакте мы получим следующее.
        
        \begin{teo}\label{teo:main theorem on compact}
            Пусть $F: Q_x \times Q_y \to \R$, $f(x) = \min_{y \in Q_y} F(x, y)$. Обозначим за $\tilde \e$ точность решения внутренней задачи \eqref{eq:inner problem}: $F(x, y_{\tilde \e}) - f(x) \le \tilde \e$. При этом, 
            \begin{itemize}
                \item $Q_x \subset \R^m$ и $Q_y \in \R^n$ -- некоторые компактные множества,
                
                \item $F(x, y)$ $\mu_x$-сильно выпуклая по $x$,
                
                \item $F(x, y)$ $\mu_y$-сильно выпуклая по $y$,
                
                \item $F(x, y)$ $L_y$-гладкая по $y$,
                
                \item $F(x, y)$ $L_{p,y}$-гладкая по $y$,
                
                \item $F(x, y)$ $L_{xy}$-гладкая по совокупности переменных функция,
                
                \item обозначим $ \delta = \frac{H}{1 - \gamma} \left( \sqrt{\frac{2}{\mu_y} D} + \sqrt{\frac{2}{\mu_y} \e} \right)^{p + 1},\ D = \max_{y \in Q_y} \left( F(x, y(x)) - F(x, y^*(x)) \right)$ и выберем $\tilde \e$ таким образом, что 
                \[
                    \delta(\tilde \e) \le \frac{\mu_x \|x_0 - x^*\|_2^2}{10 \left( \sqrt \frac{20 L_{xy}}{\mu_x} + 3 \right)}.
                \]
            \end{itemize}
            Будем решать внутреннюю задачу \eqref{eq:inner problem} проксимальным тензорным методом для сильно выпуклых функций (Теорема \ref{teo:proximal tensor method for strongly convex}). Внешнюю задачу будем решать быстрым градиентным методом с неточным оракулом для сильно выпуклых функций (Теорема \ref{teo:fgm with inexact oracle for strongly convex}). Тогда сложность конечного алгоритма составит
            \begin{equation}\label{eq:final algorithm complexity on compact simple}
                \tilde O \left( 
                    \left( \frac{L_{xy}}{\mu_x} \right)^\frac{1}{2} \left(  
                        \frac{L_{p,y} R_y^{p - 1}}{\mu_y}
                    \right)^\frac{1}{p + 1}
                \right)
            \end{equation}
            где $R_x = \|x_0 - x^*\|_2,\ R_y = \|y_0 - y^*\|_2$, а под тильду в $\tilde O$ мы занесли все логарифмические факторы.
        \end{teo}
        \begin{proof}
            Следует непосредственно из Теорем \ref{teo:proximal tensor method for strongly convex}, \ref{teo:about (delta, L) oracle and smoothness on compact}, \ref{teo:fgm with inexact oracle for strongly convex}.
        \end{proof}\qed

        \begin{rem}\label{rem:fgm vs our method}
            Сравним полученный результат с традиционным подходом в выпуклой оптимизации -- быстрым градиентным методом по совокупности переменных.
            Предположим, что $\mu_x \ge \mu_y$.
            Если решать всю задачу \eqref{eq:min-min statement last} быстрым градиентным методом, то сложность составит
            \begin{equation}
                \tilde O \left( \sqrt \frac{ L_{xy} }{ \min{\{\mu_x, \mu_y\}} } \right) = \tilde O \left( \sqrt \frac{ L_{xy} }{ \mu_y } \right).
            \end{equation}
            
            Допустим, внутренняя переменная определена на всём пространстве.
            Сравним данную оценку с \eqref{eq:final algorithm complexity simple} с точки зрения коэффициентов сильной выпуклости:
            \[
                \frac{1}{\mu_y^\frac{1}{2}} \vee \frac{1}{\mu_x^{\frac{1}{2}} \mu_y^\frac{1}{4}} \Leftrightarrow \frac{1}{\mu_y^\frac{1}{4}} \vee \frac{1}{\mu_x^{\frac{1}{2}}}.
            \]
            Таким образом, предлагаемый нами в Теореме \ref{teo:main theorem on R^n} метод выигрывает по сложности у быстрого градиентного метода, если $\mu_x \ge \sqrt{\mu_y}$.
        \end{rem}
        
        \begin{rem}\label{rem:fgm vs our method on compact}
            Теперь, допустим, внутренняя переменная определена на компакте. Сравним оценку быстрого градиентного метода с \eqref{eq:final algorithm complexity on compact simple}:
            \[
                \frac{1}{\mu_y^\frac{1}{2}} \vee \frac{1}{\mu_x^{\frac{1}{2}} \mu_y^\frac{1}{p + 1}} \Leftrightarrow \frac{1}{\mu_y^\frac{p - 1}{2(p + 1)}} \vee \frac{1}{\mu_x^\frac{1}{2}}.
            \]
            Таким образом, предлагаемый нами в Теореме \ref{teo:main theorem on compact} метод выигрывает по сложности у быстрого градиентного метода, если $\mu_x \ge \mu_y^\frac{p - 1}{p + 1}$. 
            
            К примеру, для $p = 3$ получаем $\mu_x \ge \sqrt \mu_y$ так же как и в Замечании \ref{rem:fgm vs our method}. Однако, здесь для решения внутренней подзадаче нам нужно будет вычислять производные третьего порядка в отличии от метода из Замечания \ref{rem:fgm vs our method}.
        \end{rem}
    
\paragraph{Заключение}

    В данной работе мы предложили методы для решения задач типа min-min для сильно-выпуклых по обеим переменным функций. Мы разделили задачу на внутреннюю и внешнюю подзадачи. В случае если внутренняя задача определена на всём пространстве, мы воспользовались супербыстрым тензорным методом Нестерова \cite{nesterov2021superfast}, если же она определена на компакте, то проксимальным тензорным методом \cite{ahookhosh2021high_a}. Внешняя задача решалась быстрым градиентным методом с неточным оракулом на компакте. Сложность полученного метода позволяет говорить о его преимуществах по сравнению с обычным быстрым градиентным методом в некотором специальном сеттинге (см. Замечания \ref{rem:fgm vs our method} и \ref{rem:fgm vs our method on compact}).
    
    В конце хотелось бы отметить некоторые моменты, которые не удалось проработать в рамках данной статьи. Во-первых, вместо $(\delta, L)$-оракула можно попробовать рассмотреть $(\delta, L, \mu)$-оракул. Тогда не было бы необходимости в применении рестартов к быстрому градиентному методу. Ведь, как известно, рестарты на практике дают сильно завышенную оценку количества итераций. Также, при применении $(\delta, L, \mu)$-оракула не возникает накопление ошибки (см. Теорему 5.14 в \cite{devolder2013exactness}), которое мы наблюдали в \eqref{eq:fgm with inexact oracle complexity}. 
    Во-вторых, также имеет смысл использовать методы нулевого порядка для решения внешней задачи вместо градиентных методов. В третьих, в случае, еслп внутренняя задача определена на компакте, можно для ее решения использовать супербыстрый тензорный метод (см. \cite{ahookhosh2021high_b}). В конце концов, можно рассмотреть различные обобщения постановки: равномерную выпуклость, условия Гёльдера вместо условий Липшица. Если рассматривать возможные направления исследований задачи min-min и абстрагироваться от конкретных методов, использованных внутри смешанного оракула, можно обратить внимание на отсутствие нижних оценок для данной задачи.

\bibliography{bibliography.bib}

\end{document}